\documentclass[12pt]{amsart}

\usepackage{mathrsfs}

\usepackage{amsmath,amsthm,amssymb}
\font\teneufm=eufm10
\font\seveneufm=eufm7
\font\fiveeufm=eufm5
\newfam\frakturfam

\textfont\frakturfam=\teneufm
\scriptfont\frakturfam=\seveneufm
\scriptscriptfont\frakturfam=\fiveeufm

\newtheorem{pr}{Proposition}

\newtheorem{lm}{Lemma}
\newtheorem{theor}{Theorem}
\newtheorem{co}{Corollary}

\def\bee{\begin{eqnarray}}
\def\bes{\begin{eqnarray*}}
\def\eee{\end{eqnarray}}
\def\ees{\end{eqnarray*}}

\def\b{\beta}
\def\g{\gamma}

\def\s{\sigma}
\def\t{\tau}

\def\Proof{{\sl Proof.}\ }

\title{Automorphisms of free braided associative algebras in two variables}

\begin{document}

\date{}
\maketitle

\bigskip
\begin{center}
{\bf Riza Mutalip}
\footnote {L.N. Gumilyov Eurasian National University, Nur-Sultan, Kazakhstan, e-mail: {\em mutalipriza@yahoo.com}},
{\bf Altyngul Naurazbekova}
\footnote{L.N. Gumilyov Eurasian National University, Nur-Sultan, Kazakhstan,  e-mail: {\em altyngul.82@mail.ru}},
and
{\bf Ualbai Umirbaev}\footnote{
 Wayne State University,
Detroit, MI 48202, USA, and Institute of Mathematics and Modeling, Almaty, Kazakhstan,
e-mail: {\em umirbaev@wayne.edu}}
\end{center}

\begin{abstract} We describe the groups of automorphisms of two generated free braided associative algebras
with involutive diagonal braidings over a field of characteristic $\neq 2$. Depending on the form of the diagonal involutive braiding, five different automorphism groups arise as automorphism groups of two generated free braided associative algebras. This list covers all three automorphism groups that arise in the case of quantum planes \cite{AC}.
\end{abstract}

\noindent
{\bf Mathematics Subject Classification (2010):} 16W20, 16T20, 16T25.

\noindent
{\bf Key words:} Yang-Baxter equation, braiding, free associative algebra,  automorphism.

\section{Introduction}

\hspace*{\parindent}

It is well known \cite{Czer,Jung,Kulk,Makar} that all automorphisms of the polynomial algebra $K[x_1,x_2]$ and the
free associative algebra $K\langle x_1,x_2\rangle$ in two variables over an arbitrary field $K$ are tame. Moreover
\cite{Czer,Makar}, the automorphism groups of the algebras $K[x_1,x_2]$ and $K\langle x_1,x_2\rangle$ are isomorphic,
that is
\bes
\mathrm{Aut} K[x_1,x_2]\cong \mathrm{Aut} K\langle x_1, x_2\rangle.
\ees

It is also known that all automorphisms of two generated free Poisson algebras over a field of characteristic zero \cite{MLTU}
and two generated free right-symmetric algebras over an arbitrary field \cite{KMLU} are tame.
P. Cohn \cite{Cohn} proved that all automorphisms of free Lie algebras of a finite rank over a field are tame.
An analogue of this theorem is true for free algebras of any homogeneous Schreier variety of algebras \cite{Lewin}.
Recall that the varieties of all non-associative algebras \cite{Kurosh}, commutative and anti-commutative algebras \cite{Shir},  Lie algebras  \cite{Shir1, Witt} and  Lie superalgebras \cite{Mikh, Shtern} are Schreier.

Let $q=(q_{ij})\in M_n(K)$ be a matrix of order  $n\geq 2$ whose elements $q_{ij}\in K^*$ satisfy the relations
$q_{ii}=q_{ij}q_{ji}=1$  \text{for all} $1\leq i,j\leq n$. Denote by $ \mathrm{O}_q=K_q[x_1,\ldots,x_n] $
the associative algebra with unity generated by the elements $x_1,\ldots,x_n$ which satisfy the relations
\bes
x_ix_j=q_{ij}x_jx_i, \ \  \  \ 1\leq i, j \leq n.
\ees
The algebra $O_q$ is called a {\em quantum polynomial algebra} \cite {Artamonov02}. In the case $n=2$ the algebra
$O_q$ is called a {\em quantum plane}.

The automorphism groups of the quantum polynomial algebras were described by J. Alev and M. Chamarie \cite{AC} (see also \cite {Artamonov02}).
All automorphisms of the quantum plane $O_q$ with $q_{12}\neq\pm1$ are
 {\em toric}, i.e. are isomorphic to the multiplicative group $K^*\times K^*$. Moreover,
the automorphism  group $\mathrm{Aut}(O_q)$ of the quantum plane $O_q$ is isomorphic to
$\mathrm{Aut} K\langle x_1,x_2 \rangle$ if $q_{12}=1$ and  is isomorphic to the semidirect product $(K^*\times K^*) \rtimes \mathbb{Z}_2$ if $q_{12}=-1$.
Automorphisms of the field of fractions of quantum polynomial algebras were studied in \cite{AD,Artamonov00,Artamonov99}.

A {\it braided space} is a linear space over an arbitrary field $K$ with a linear
map $\tau :V\otimes V\rightarrow V\otimes V$
that satisfies the {Yang--Baxter equation} or the {\it braid relation}
\begin{equation} \label{f2}
(\tau\otimes {\rm id})({\rm id}\otimes \tau )(\tau\otimes {\rm id})
=({\rm id}\otimes \tau )(\tau\otimes {\rm id})({\rm id}\otimes \tau ).
\end{equation}
The linear map $\tau $ will be called a {\it braiding} of the space $V$.

If $ x_1, x_2, \ldots, x_n $ is a basis of a linear space $V$, then for arbitrary parameters
$ q_ {is} \in K $,  $ 1 \leq i, s \leq n $,  the linear map $\tau$ defined by
\begin{equation}\label{f3}
\tau: x_i\otimes x_s\mapsto q_{is}\cdot  x_s\otimes x_i
\end{equation}
is a braiding and is called a {\em diagonal braiding}. Denote by $\tau=(q_{ij})\in M_n(K)$ this diagonal braiding.

A braiding $\tau :V\otimes V\rightarrow V\otimes V$ is called {\it involutive} if $\t^2=\mathrm{id}$. The diagonal braiding $\tau=(q_{ij})\in M_n(K)$
is involutive if and only if $q_{ij}q_{ji}=1$  for all $1\leq i,j\leq n$. In particular, $q_{ii}=\pm 1$ for all $i$. The best known example of an involutive braiding is the {\em ordinary flip} $\theta$ defined by
$\theta(x\otimes y)=y\otimes x$ for all $x,y \in V$.

It is well known  \cite{Khar15} that every braiding  $\tau : V\otimes  V\rightarrow V\otimes  V$ of a vector space $V$ can be uniquely extended to a braiding $\tau'$ of the tensor algebra
$$T(V)=\bigoplus_{i=0}^{\infty } V^{\otimes i}.$$
Moreover,
there is a canonical braided structure of the Hopf algebra on $T(V)$, which plays an important
role in quantum Lie theory (see, for example \cite{Khar15,Tak00}).  Analogues of these results for free nonassociative algebras are also true \cite{Umi34}.
If $x_1,x_2,\ldots,x_n$ is a basis of the linear space $V$ then $T(V)$ is isomorphic to the free associative algebra $K\langle x_1,x_2,\ldots,x_n \rangle$ over $K$ freely generated by $x_1,x_2,\ldots,x_n$.

In this paper we describe the groups of automorphisms of free braided associative algebras in two variables $A_{\tau}=(K\langle x_1,x_2\rangle, \tau)$  with a diagonal involutive braiding $\tau = (q_{11},q_{12},q_{21},q_{22})$ over a field of characteristic $\neq 2$. If $q_{ij}=1$ for all $i,j$ then $\tau=\theta$ is the ordinary flip and the group of automorphisms $\mathrm{Aut}A_{\tau}$ of $A_{\tau}$ coincides with the group of all automorphisms $\mathrm{Aut}K\langle x_1,x_2\rangle$ of the free associative algebra $K\langle x_1,x_2\rangle$. If $q_{ij}=-1$ for all $i,j$ then the group of automorphisms $\mathrm{Aut}A_{\tau}$ coincides with the group of all odd automorphisms $G_{\mathrm{odd}}$ (full definitions are given in Section 3) of the free associative algebra $K\langle x_1,x_2\rangle$. Furthermore, $\mathrm{Aut}A_{\tau}$ is isomorphic to the toric group $K^*\times K^*$ if $q_{12}\neq \pm 1$ or $q_{12}=-1, q_{11}q_{22}=-1$, and isomorphic to the group $(K^*\times K^*) \rtimes \mathbb{Z}_2$ if $q_{11}=q_{22}, q_{12}=q_{21}, q_{11}q_{12}=-1$. If $q_{12}=1$ and $q_{11}q_{22}=-1$ then $\mathrm{Aut}A_{\tau}$ is isomorphic to the subgroup of triangular automorphisms of $K\langle x_1,x_2\rangle$ of the form
\bes
\{\varphi\in \mathrm{Aut}K\langle x_1,x_2\rangle \mid \varphi(x_1)=\alpha x_1+g(x_2^2), \varphi(x_2)=\beta x_2, \alpha,\beta\in K^*, g(x)\in K[x]\}.
\ees

The paper is organized as follows. In Section 2 some necessary definitions and facts on free braided associative algebras are given.
In Section 3 we classify involutive diagonal braidings on the free associative algebra $K\langle x_1,x_2\rangle$ up to isomorphism.
In Section 4 we prove the main result of the paper on the description of all groups of automorphisms of the two generated free braided associative algebras
over a field of characteristic $\neq2$
with an involutive diagonal braiding .

\section{Braided algebras}
\hspace*{\parindent}

A {\em braid monoid }  $B_n$ \cite {Khar15}  is an associative monoid generated by braids $s_1, s_2, \ldots, s_{n-1}$
subject to the relations
\begin{equation}\label{f1}
s_ks_{k+1}s_k=s_{k+1}s_ks_{k+1}, \, \, s_is_j=s_j s_i, \,\, 1\leq k<n-1, \, \, |i-j|> 1.
\end{equation}

The group with generators  $ s_1, s_2, \cdots,  s_ {n-1} $ and defining relations
 (\ref{f1}) is called  the {\em Artin braid group}.

Let $V$ be a linear space with a braiding $ \tau: V \otimes V \rightarrow V \otimes V $.
Consider the linear maps
$$
\tau _i={\rm id}^{\otimes (i-1)} \otimes \tau \otimes {\rm id}^{\otimes (n-i-1)}:
V^{\otimes n}\rightarrow V^{\otimes n}, \ \ 1\leq i<n.
$$

By (\ref {f2}), the maps $\tau_i$ satisfy all defining relations (\ref{f1}) of the braid monoid, i.e.,
$$\tau_i\tau_{i+1}\tau_i=\tau_{i+1}\tau_{i}\tau_{i+1}, \ \  1\leq i<n-1; \ \ \tau_i\tau_j=\tau_j\tau_i, \ \  |i-j|>1.$$

Put $[k;k]=1$ and
$$ [m;k]=\tau_{k-1}\tau_{k-2}\cdots \tau_{m+1}\tau_m, \ \ [k;m]=\tau_m\tau_{m+1}\cdots \tau_{k-2}\tau_{k-1}, \ \ m<k.$$

Consider a map $\nu_r^{k, n}: V^{\otimes n}\rightarrow V^{\otimes n}, \, k \leq r < n$, defined as a superposition of the  $\tau_i$'s:
\begin{equation}\label{f4}
\nu_r^{k,n}=[k;r+1] [k+1;r+2] \cdots [k+n-r-1;n].
\end{equation}
It is proven in \cite{Khar15}, that
$$\nu _r^{k,n}=[n;r][n-1;r-1]\cdots [n-r+k;k].$$

Let $V$ and $V^{\prime} $ be spaces with braidings $ \tau $ and $ \tau^{\prime}$, respectively.
A linear map $ \varphi: V \rightarrow V^{\prime} $ is called a {\it homomorphism of braided spaces} if
$$
\tau (\varphi \otimes \varphi )=(\varphi \otimes \varphi )\tau^{\prime} .
$$

An algebra $R$ with a multiplication ${\bf m}: R \otimes R \rightarrow R $ is called a {\em braided algebra}
if $R$ is a braided space and
$$({\bf m}\otimes {\rm id})\tau =\tau_2\tau_1({\rm id}\otimes{\bf m}), \, \, \,
 ({\rm id}\otimes {\bf m})\tau =\tau_1\tau_2({\bf m}\otimes {\rm id}).$$

As above, in these formulas we use the so-called ``exponential notation'' for actions of the operators,
that is, the operators in a superposition act from the left to the right.

A {\em homomorphism of braided algebras} is a linear map that is both a homomorphism of algebras and braided spaces.

Let $V$ be a linear space and
$$T(V) = \bigoplus_{i=0} ^{\infty} V ^{\otimes i}$$
be its tensor algebra.
The product in the tensor algebra $T(V)$ will be denoted by ${\bf m}$, that is
\bes
{\bf m} : T(V)\otimes' T(V)\longrightarrow T(V),
\ees
where the sign $\otimes'$ is the same tensor product $\otimes$ with one additional function that separates the two tensors to which the product  ${\bf m}$ is applied. For example,
 $(u\otimes' v){\bf m}=u\otimes v$. Very often we write $v_1v_2\ldots v_k$ instead of $v_1\otimes v_2\otimes \ldots \otimes v_k\in V^{\otimes k}$.

We fix a linear basis $X = \{x_i | i \in I \}$ of $V$. Then $T(V)=K\langle X \rangle$
is the free associative algebra over $K$ with a free set of generators $X$. The set of all associative words $X^*$ in the alphabet $X$ forms a linear basis for $K\langle X \rangle$. If $w=x_{i_1}x_{i_2}\ldots x_{i_k}\in X^*$ then the length $k$ of $w$ will be denoted by $\mathrm{d}(w)$. Notice that $w\in V^{\otimes k}$. Set $\mathrm{mdeg}(x_i)=\varepsilon_i$, where $\varepsilon_1, \ldots, \varepsilon_n$ is the standard basis for $Z^n$. The {\it multidegree} of the word $w$ is defined by
\bes
\mathrm{mdeg}(v)=\mathrm{mdeg}(x_{i_1})+\mathrm{mdeg}(x_{i_2})+\ldots+\mathrm{mdeg}(x_{i_s}).
\ees

It is well known  \cite{Khar15} that every braiding  $\tau$ has a unique extension $ \tau'$ on the free associative algebra
 $K \langle X \rangle $ such that  $K \langle X \rangle $ is a braided algebra. For any $0 \leq r\leq n$ denote by
 $\theta_r$ the linear map

$$\theta_r : V^{\otimes n}\rightarrow V^{\otimes r} \otimes' V^{\otimes (n-r)}$$
defined by
$$
(z_1z_2\ldots z_n)\theta_r=z_1z_2\ldots z_r\otimes' z_{r+1}\ldots z_{n}, \ \ \ z_i\in X.
$$

The braiding $\tau'$ is defined \cite{Khar15} by
\begin{equation} \label{f5}
(u\otimes' v)\tau^{\prime } =(u\otimes v) \nu_r^{1,n} \, \theta _{n-r}, \ \ \ u\in V^{\otimes r}, \ v\in V^{\otimes (n-r)}.
\end{equation}

In particular,
$ (1\otimes' v)\tau^{\prime }=v\otimes' 1$ and
$(u\otimes' 1)\tau^{\prime }=1\otimes' u$ by this definition.

\begin{lm}\label{1}
Let $V$ be a vector space over a field $K$ with a linear basis $X =\{x_1, x_2, \ldots, x_n\}$ and let $\tau$ be a diagonal braiding of $V$ defined by (\ref{f3}). Let  $u,v\in X^*$, $\mathrm{mdeg}(u)=(s_1,s_2,\ldots,s_n)$, and $\mathrm{mdeg}(v)=(t_1,t_2,\ldots,t_n)$.
Then
$$ (u\otimes' v)\tau'=\prod_{i,j} q_{ij}^{s_it_j}(v\otimes' u).$$
\end{lm}
\Proof
Set $s=s_1+\cdots +s_n=\mathrm{d}(u)$ and $t=t_1+\ldots +t_n=\mathrm{d}(v)$.
We prove the statement of the lemma by induction on $t$.
If $t=0$, then $v=1$ and  $(u\otimes' 1)\tau'=1\otimes' u$. Suppose that $t\geq 1$ and $v=v'x_{h}$. Set $\mathrm{mdeg} (v')=(l_1, l_2, \cdots, l_n)=(t_1, t_2, \ldots, t_{h-1}, t_h-1, t_{h+1}, \ldots, t_n)$.
 By (\ref{f5}) and (\ref{f4}), we have
\bes
(u\otimes' v)\tau'=(u\otimes' v'x_h)\tau'=(u\otimes v'\otimes x_h)\nu^{1, s+t}_s\theta_t\\
=(u\otimes v'\otimes x_h)[1, s+1][2, s+2]\ldots [t-1, s+t-1][t, s+t]\theta_t\\
=(u\otimes v'\otimes x_h)\nu^{1, s+t-1}_s (\tau_{s+t-1}\tau_{s+t-2}\cdots\tau_t)\theta_t\\
=(((u\otimes v') \nu^{1, s+t-1}_s \theta_{t-1})\otimes x_h)(\tau_{s+t-1}\tau_{s+t-2}\cdots\tau_t)\theta_t\\
=(((u\otimes' v')\tau')\otimes x_h)(\tau_{s+t-1}\tau_{s+t-2}\cdots\tau_t)\theta_t.
\ees

By the induction proposition, we get
\bes
(u\otimes' v)\tau'= (u\otimes' v'\otimes x_h)\tau'
=\prod_{i,j} q_{ij}^{s_il_j}( v'\otimes u\otimes x_h)
(\tau_{s+t-1}\tau_{s+t-2}\cdots\tau_t)\theta_t\\
=\prod_{i,j} q_{ij}^{s_il_j}q_{ih}^{s_i}(v'\otimes x_h\otimes u)\theta_t
=\prod_{i,j} q_{ij}^{s_il_j}q_{ih}^{s_i}(v\otimes u)\theta_t=
\prod_{i,j} q_{ij}^{s_it_j}(v\otimes' u). \ \ \ \Box
\ees

\section{Diagonal braidings on $K\langle x_1,x_2\rangle$}
\hspace*{\parindent}

Let $V$ be a two-dimensional vector space over a field $K$ with a linear basis $x_1,x_2$. Denote by
\bes
\tau = (q_{11},q_{12},q_{21},q_{22})
\ees
 the diagonal braiding of $V$ defined by (\ref{f3}).

  Let $A=K\langle x_1, x_2 \rangle=T(V)$ be the free associative
algebra of the vector space $V$.  Then $A$  is a free braided associative algebra with the braiding  $\tau'$ defined in Section 2. Denote this extension by the same symbol $\tau$ and denote by $A_{\tau} = (A, \tau)$ the free associative algebra $A$ with braiding $\tau$.

By interchanging the variables $x_1$ and $x_2$, from $\tau$ we can get another diagonal braiding, which will be denoted by
\bes
\tau^* = (q_{22},q_{21},q_{12},q_{11}).
\ees

\begin{pr}\label{p1}
Let  $A_{\tau}=(K\langle x_1,\, x_2\rangle, \tau)$ and
 $A_{\sigma}=(K\langle x_1,\, x_2\rangle, \sigma)$ be  free braided associative algebras
in two generators $x_1,x_2$ over an arbitrary field $K$ with diagonal braidings
 $\tau$ and $\sigma$, respectively. Then
$ A_\tau \cong  A _\sigma$ if and only if $\sigma = \tau$ or $ \sigma = \tau^*$.
\end{pr}
\Proof
Let $ \varphi: A_\tau \rightarrow A_\sigma$  be an isomorphism  and $\varphi(x_1)=f_1$,  $\varphi(x_2)=f_2$. Denote by
 $L(f_1)=\alpha_1x_1+ \alpha _2 x_2$ and  $L(f_2)= \beta_1 x_1+ \beta_2 x_2$  the linear parts of $f_1$ and $f_2$,  respectively. Since $\varphi$ is an automorphism, we have
\begin{equation}\label{f6}
det\left[\begin{array}{cccc}
\alpha_1 & \alpha_2 \\
\beta_1 & \beta_2 \\
\end{array}\right]\neq 0
\end{equation}
and
$${((x_i\otimes x_j)\tau)(\varphi \otimes \varphi)}=((x_i\otimes x_j)(\varphi \otimes \varphi))\sigma,
  \, 1\leq i, j\leq 2.$$

Let $\tau = (q_{11},q_{12},q_{21},q_{22})$ and $\sigma = (p_{11},p_{12},p_{21},p_{22})$. Then
\begin{equation}\label{f7}
q_{ij}(\varphi(x_j)\otimes\varphi(x_i))=(\varphi(x_i)\otimes\varphi(x_j))\s, \, \, 1 \leq i,  j \leq 2.
\end{equation}
Consequently, for $i=1$ and $j=2$ we have
$$L(q_{12} (\varphi(x_2)\otimes \varphi(x_1))) =L((\varphi (x_1) \otimes \varphi (x_2))\s),$$
that is
$$q_{12}(\beta_1x_1+\beta_2x_2)\otimes(\alpha_1x_1+\alpha_2x_2)= ((\alpha_1x_1+\alpha_2 x_2)\otimes(\beta_1x_1+\beta_2x_2))\s.$$
Hence
$$q_{12}(\alpha_1\beta_1(x_1\otimes x_1)+\alpha_2\beta_1(x_1\otimes x_2)+\alpha_1\beta_2(x_2\otimes x_1)+\alpha_2\beta_2(x_2\otimes x_2))$$
$$=p_{11}\alpha_1\beta_1(x_1\otimes x_1)+p_{12}\alpha_1\beta_2(x_2\otimes x_1)+p_{21}\alpha_2\beta_1(x_1\otimes x_2)+p_{22}\alpha_2\beta_2(x_2\otimes x_2).$$
Comparing the coefficients of the terms $x_1\otimes x_1, \,x_1\otimes x_2, \,   x_2\otimes x_1,\,  x_2\otimes x_2$
 in this equality, we get
\begin{equation}\label{f8}
(p_{11}-q_{12})\alpha_1\beta_1=(p_{21}-q_{12})\alpha_2\beta_1=(p_{12}-q_{12})\alpha_1\beta_2=(p_{22}-q_{12})\alpha_2\beta_2=0.
\end{equation}
Applying (\ref{f7}) for the remaining values of $i$ and $j$, we similarly obtain that
\begin{equation}\label{f9}
(p_{11}-q_{11})\alpha^2_1=(p_{21}-q_{11})\alpha_1\alpha_2=(p_{12}-q_{11})\alpha_1\alpha_2=(p_{22}-q_{11})\alpha^2_2=0,
\end{equation}
\begin{equation}\label{f10}
(p_{11}-q_{21})\alpha_1\beta_1=(p_{21}-q_{21})\alpha_1\beta_2=(p_{12}-q_{21})\alpha_2\beta_1=(p_{22}-q_{21})\alpha_2\beta_2=0,
\end{equation}
 \begin{equation}\label{f11}
(p_{11}-q_{22})\beta^2_1=(p_{21}-q_{22})\beta_1\beta_2=(p_{12}-q_{22})\beta_1\beta_2=(p_{22}-q_{22})\beta^2_2=0.
\end{equation}

By (\ref{f6})  $\alpha_1 \beta_2 \neq 0$ or $\alpha_2\beta_1\neq0$. If $\alpha_1 \beta_2 \neq 0$ then it follows from
(\ref{f8}) -- (\ref{f11}) that
$$ p_{12}=q_{12},\, p_{11}=q_{11}, \, p_{21}=q_{21}, \, p_{22}=q_{22}.$$
If $\alpha_2\beta_1\neq0$ then, similarly, we get
 $$ p_{21}=q_{12}, \, p_{22}=q_{11}, \, p_{12}=q_{21}, \, p_{11}=q_{22}.$$
Consequently, if  $ A_\tau\cong A_\sigma$, then $\sigma=\tau$ or $\sigma=\tau^*$.

Let $\sigma = \tau^*$.  Consider the automorphism $\varphi$ of the free associative algebra $K\langle x_1,\, x_2\rangle$ given by
 $\varphi(x_1)= x_2$ and $\varphi(x_2)= x_1$.
 To prove that $\varphi: A_\tau \rightarrow A_\sigma$ is an isomorphism of braided algebras, it is sufficient to check that
$\tau(\varphi\otimes \varphi)=(\varphi\otimes\varphi)\tau^* $.
Let $ u, v \in X^*, \,   s_i=\deg_{x_i}u, \, t_i=\deg_{x_i}v , \,  i\in \{1,2\}$. Note that
$ \deg_{x_1}\,\varphi(u)=s_2, \, \deg_{x_2}\,\varphi(u)=s_1, \, \deg_{x_1}\,\varphi(v)=t_2, \, \deg_{x_2}\,\varphi(v)= t_1$.
 By Lemma 1, we have
$$((u\otimes v)\tau)(\varphi\otimes\varphi)=q_{11}^{s_1t_1}q_{12}^{s_1t_2}q_{21}^{s_2t_1}q_{22}^{s_2t_2}(\varphi(v)\otimes \varphi(u))$$ and
$$((u\otimes v )(\varphi\otimes\varphi))\tau^*=q_{11}^{s_1t_1}q_{12}^{s_1t_2}q_{21}^{s_2t_1}q_{22}^{s_2t_2}(\varphi(v)\otimes \varphi(u)).$$
Hence
$$((u\otimes v)\tau)(\varphi\otimes\varphi)=((u\otimes v )(\varphi\otimes\varphi))\tau^*$$
and  $A_\tau\cong A_\s$. $\Box$

Notice that the braiding $\tau$ from (\ref{f3}) is involutive if and only if
\bes
q_{ij}q_{ji}=1
\ees
for all $1\leq i,j\leq n$. Consequently, in the case of two variables the braiding $\tau$ from (\ref{f6}) is involutive if and only if
$$ q_{11}^2=q_{12}q_{21}=q_{22}^2=1.$$
In particular, $q_{11}=\pm 1$ and $q_{22}=\pm 1$.
 Note that if $\tau$ is involutive then $\tau^*$ is also involutive.

\begin{co}\label{c1}
Let  $A_{\tau}=(K\langle x_1,\, x_2\rangle, \tau)$ be a free braided associative algebra in two variables $ x_1,x_2 $ over a field $K$ of characteristic  $p\neq2$ with an involutive diagonal braiding
 $\tau=(q_{11},q_{12},q_{21},q_{22})$.
Then either $q_{12}\neq \pm 1$ or $A_\tau $ is isomorphic to  $A_\s$, where $\s$ is equal to one of the following vectors:
\bes
(1,1,1,1), (-1,-1,-1,-1),  (1,1,1,-1),   (-1,1,1,-1), (1,-1,-1,1), (1,-1,-1,-1).
\ees
\end{co}
\Proof
If $q_{12}=1$ then $\tau$ is involutive if only if it is equal to one of
 the following vectors:
\bes
(1,1,1,1), (1,1,1,-1), (-1,1,1,1), (-1,1,1,-1).
\ees
  By Proposition \ref{p1},  $A_{(-1,1,1,1)}\cong A_{(1,1,1,-1)}$.

If $q_{12}=-1$  then  $\tau$ is equal to one of the following vectors:
\bes
(1,-1,-1,1), (1,-1,-1,-1), (-1,-1,-1,1), (-1,-1,-1,-1).
\ees
By Proposition \ref{p1}, $A_{(1,-1,-1,-1)}\cong A_{(-1,-1,-1,1)}$. $\Box$

\section{Automorphisms of  $A_{\tau}=(K\langle x_1,\, x_2\rangle, \tau)$}

\hspace*{\parindent}

Let  $K\langle X \rangle = K\langle x_1,x_2,\ldots,x_n \rangle$  be the free associative algebra over $K$
with the set of free generators $X=\{x_1,x_2,\ldots,x_n\}$.
 Denote by $\varphi=(f_1,f_2,\ldots,f_n)$ an automorphism of the algebra  $K\langle X \rangle$ such that
$ \varphi(x_1)=f_1,\varphi(x_2)=f_2,\ldots,\varphi(x_n)=f_n$. 
Denote by $\mathrm {Aut} K \langle X \rangle $ the set of all automorphisms of the algebra $K\langle X\rangle$ over $K$.
An automorphism of the form
$$(x_1,\ldots,x_{i-1},\alpha x_i+f, x_{i+1},\ldots,x_n),$$
where  $0\neq\alpha\in K,\ f\in K \langle X\setminus \{x_i\}\rangle$, is called {\em  elementary}. An automorphism $\varphi$ is called {\em tame}
if it is a product of elementary automorphisms.

Consider the natural grading
 $$K\langle X\rangle = A_0 \oplus A_1 \oplus \cdots \oplus A_m \oplus \cdots, $$
where $A_0 = K$, $A_1 = Kx_1 + Kx_2+\cdots +Kx_n$, and $A_m$ is the linear span of words of length $m$.
Every nonzero element $f\in K \langle X \rangle$ can be uniquely represented in the form
$$f=f_0+f_1+\cdots + f_{m-1}+f_m,$$
 where  $f_i \in A_i$,  $f_m\neq 0$. Then $f_m$ is called {\it the highest homogeneous part} of  $f$ and is denoted by
$\overline{f}$. Set also $\deg f=m$. The degree of an automorphism $\varphi=(f_1,f_2,\ldots,f_n)$ is defined by
$$\deg\varphi=\deg f_1+\deg f_2+\ldots+\deg f_n.$$

An {\em elementary transformation} of an $n$-tuple $\varphi=(f_1,f_2,\ldots,f_n)$ changes only one element $f_i$ for some $i$ with an element of the form $\alpha f_i+g$,
where $0\neq \alpha \in K$ and $g \in K \langle  f_j | j \neq i \rangle$, and leaves unchanged all other elements $f_j$, where $j\neq i$. We write $\varphi \rightarrow \psi$ if $\psi$ is obtained from $\varphi$ by one elementary transformation.
An automorphism  $\varphi$ is called  {\em elementarily reducible} if there exists an automorphism $\psi$ such that
$\varphi \rightarrow \psi$ and $\deg\psi < \deg\varphi$.

Every automorphism of the free associative algebra $A=K\left\langle  x_1, x_2\right\rangle$ in two variables over a field $K$ is tame
\cite{Czer, Makar}.
In the next lemma we formulate a property of automorphisms of $A=K\left\langle  x_1, x_2\right\rangle$ that we need in the future.
\begin{lm}\label{2}
Let $A = K \langle x_1, x_2 \rangle$ be the free associative algebra in two variables $x_1,x_2$.
If $\varphi=(f_1, f_2)$ is an automorphism of the algebra $A$ with degree $\geq 3$, then
$$\bar{f_1}=\alpha(ax_1+bx_2)^r,\ \  \bar{f_2}=\beta(ax_1+bx_2)^s,$$
where $ r+s\geq 3$, $r\mid s$ or $s\mid r$, $\alpha, \beta\in K^*$, $a,b\in K$,  and $(a,b)\neq (0,0)$.
\end{lm}
\Proof
Let $\deg f_1=r$ and $\deg f_2=s$. The property $r\mid s$ or $s\mid r$ is explicitly formulated in \cite{Czer, Makar}.
Without loss of generality, we can assume that  $r\leq s$.
Prove the remaining statement of the lemma by induction on $\deg\varphi=r+s$.

Every automorphism of $A$ with degree $\geq 3$ is elementarily reducible \cite{Czer, Makar}.
Therefore, there is an elementary transformation $(f_1,f_2)\rightarrow (f_1, \g f_2-g(f_1))=(f_1,f_2')$,
where $\g \neq 0$, $g(f_1)\in K\langle  f_1\rangle$, such that $\deg f_2'=\deg(\g f_2-g(f_1)) <s$. If $\deg(f_1,f_2')\geq 3$ then by the induction proposition we get $\overline{f_1}=\alpha(ax_1+bx_2)^r$. Notice that this is true even if $\deg(f_1,f_2')=2$. Consequently,
\bes
\overline{f_2}=1/\g \overline{g(f_1)}= \b(ax_1+bx_2)^s
\ees
for some $\b\in K^*$. $\Box$

We fix the  grading
\bee\label{u1}
A= B_0\oplus B_1
\eee
of the algebra $A$, where $B_0$ is the linear span of all monomials of even length and $B_1$ is the linear span of all monomials of odd length.

An automorphism $\varphi=(f_1, f_2)$ of $A$ is called {\it odd}
if $f_1,f_2 \in B_1$. The set of all odd automorphisms forms a group and this group will be denoted by $G_{\mathrm{odd}}$.

Define the following subgroups of the group of automorphisms
$A=\mathrm {Aut}K\langle x_1,x_2 \rangle$:

(1) $G_1=\mathrm {Aut}K\langle x_1,x_2\rangle$ is the group of all automorphisms of $A$;

(2) $G_2=G_{odd}$ is the group of all odd automorphisms of $A$;

(3) $G_3=\{\varphi\in G_1 \mid \varphi=(\alpha_1x_1,\beta_2x_2) \,\text{or} \,\varphi=(\beta_1x_2, \ \alpha_2x_1), \ \alpha_1, \beta_2, \alpha_2,\beta_1 \in K^*\}$;

(4) $G_4=\{\varphi\in G_1 \mid \varphi=(\alpha_1x_1+g(x_2^2), \beta_2x_2), \, \alpha_1,\, \beta_2\in K^*, g(x)\in K[x]\}$;

(5) $G_5 =\{\varphi\in G_1 \mid \varphi =(\alpha_1x_1,\beta_2x_2), \ \alpha_1, \ \beta_2\in K^* \}$ is the group of all {\em toric} automorphisms of
 $K\langle x_1,x_2\rangle$.

Let $A_{\tau}=(A, \tau)$ be the free braided associative algebra in two generators
 $x_1,x_2$  over a field $K$ of arbitrary characteristic $\neq 2$ with an involutive diagonal braiding $\tau=(q_{11},q_{12},q_{21},q_{22})$. Notice that an automorphism
$\varphi\in \mathrm{Aut}A$ of the algebra $A$ is an automorphism of $A_{\tau}$ if and only if
\bee\label{f12}
\tau(\varphi\otimes \varphi)=(\varphi\otimes \varphi)\tau.
\eee
We identify the group of automorphisms $\mathrm{Aut}A_{\tau}$ of $A_{\tau}$ with the corresponding subgroup of $\mathrm{Aut}K\langle x_1,x_2\rangle$.

\begin{lm}\label{3}
 If $\tau=(1,1,1,1)$ then $\mathrm {Aut} A_\tau = G_1$.
\end{lm}
\Proof For any $u,\ v \in X^*$ we have $(u\otimes v)\tau=v\otimes u$ by Lemma \ref{1}.
Then $$((u\otimes v)\tau)(\varphi\otimes \varphi)=\varphi(v)\otimes \varphi(u)=(\varphi(u)\otimes \varphi(v))\tau=((u\otimes v)(\varphi\otimes \varphi))\tau$$
for any $\varphi\in G_1$ and any $u,\,v \in X^*$.
Consequently, (\ref{f12}) is true for all $\varphi\in G_1$. This means that $G_{\tau} = G_1$.  $\Box$

\begin{lm}\label{4}
If $\tau=(-1,1,1,-1)$ or $\tau=(1,-1,-1,1)$ then $\mathrm {Aut}A_\tau = G_3$.
\end{lm}
\Proof
Let $\varphi=(f_1,f_2)$ be an automorphism of the algebra $A_\tau$ with degree $\geq 3$. By Lemma 2, we have
\begin{equation}\label{f13}
\overline{f_1}=\gamma_1(ax_1+bx_2)^{s_1}, \\ \overline{f_2}=\gamma_2(ax_1+bx_2)^{s_2},
\end{equation}
where $s_1+s_2\geq3,\,s_1\mid s_2$, or $s_2\mid s_1,\,\gamma_1,\,\gamma_2\in K^*,\,a,\,b\in K$ and $(a, b)\neq 0$. By (\ref{f12}), we have
$$ \overline{((x_i\otimes x_j)\tau)(\varphi\otimes\varphi)}=\overline{((x_i\otimes x_j)(\varphi \otimes \varphi))\tau},
\, i, j \in\{1,\,2\}.$$
Hence
$$\overline{q_{ij}\varphi (x_j)\otimes \varphi (x_i)}=\overline{(\varphi(x_i)\otimes \varphi(x_j))\tau}, \, i, j \in\{1,\,2\}.$$
By (\ref{f13}), we have
$$q_{ij}\gamma_i\gamma_j(ax_1+bx_2)^{s_j}\otimes (ax_1+bx_2)^{s_i}=\gamma_i\gamma_j\left((ax_1+bx_2)^{s_i}\otimes (ax_1+bx_2)^{s_j}\right)\tau, \, i, j \in\{1,\,2\}.$$
Using Lemma 1 and comparing the coefficients of the terms $x_1^{s_j}\otimes x_1^{s_i},\,x_2^{s_j}\otimes x_1^{s_i},\,x_1^{s_j}\otimes x_2^{s_i},\,x_2^{s_j}\otimes x_2^{s_i}$ in this equality, we obtain
$$ (q_{ij}-q^{s_is_j}_{11})a^{s_i+s_j}=(q_{ij}-q^{s_is_j}_{12})a^{s_i}b^{s_j}=(q_{ij}-q^{s_is_j}_{21})a^{s_j}b^{s_i}=(q_{ij}-q^{s_is_j}_{22})b^{s_i+s_j}=0, \,\,i, j\in \{1,\,2\}. $$
Varying the values of $i, j \in\{1,\,2\}$, we get
\begin{equation}\label{f14}
(q_{11}-q^{s_1^2}_{11})a^{2s_1}=(q_{11}-q^{s_1^2}_{12})a^{s_1}b^{s_1}=(q_{11}-q^{s_1^2}_{21})a^{s_1}b^{s_1}=(q_{11}-q^{s_1^2}_{22})b^{2s_1}=0,
\end{equation}
\begin{equation}\label{f15}
(q_{12}-q^{s_1s_2}_{11})a^{s_1+s_2}=(q_{12}-q^{s_1s_2}_{12})a^{s_1}b^{s_2}=(q_{12}-q^{s_1s_2}_{21})a^{s_2}b^{s_1} =(q_{12}-q^{s_1s_2}_{22})b^{s_1+s_2}=0,
\end{equation}
\begin{equation}\label{f16}
(q_{21}-q^{s_1s_2}_{11})a^{s_1+s_2}=(q_{21}-q^{s_1s_2}_{12})a^{s_2}b^{s_1} =(q_{21}-q^{s_1s_2}_{21})a^{s_1}b^{s_2}=(q_{21}-q^{s_1s_2}_{22})b^{s_1+s_2}=0,
\end{equation}
\begin{equation}\label{f17}
(q_{22}-q^{s_2^2}_{11})a^{2s_2}=(q_{22}-q^{s_2^2}_{12})a^{s_2}b^{s_2}=(q_{22}-q^{s_2^2}_{21})a^{s_2}b^{s_2}=(q_{22}-q^{s_2^2}_{22})b^{2s_2}=0.
\end{equation}

If $\tau=(-1, 1, 1, -1)$ then (\ref{f14}) -- (\ref{f17}) implies that
\bes
(-1-(-1)^{s^2_1})a^{2s_1}=(-1-(-1)^{s^2_1})b^{2s_1}\\
=(1-(-1)^{s_1s_2})a^{s_1+s_2}=(1-(-1)^{s_1s_2})b^{s_1+s_2}\\
=(-1-(-1)^{s^2_2})a^{2s_2}=(-1-(-1)^{s^2_2})b^{2s_2}=0.
\ees
It is easy to check that over a field of characteristic $\neq 2$, these equations give $a=b=0$.

If $\tau=(1, -1, -1, 1)$ then  (\ref{f15}) gives
 $$(-1-1^{s_1s_2})a^{s_1+s_2}=(-1-1^{s_1s_2})b^{s_1+s_2}=0$$
and, consequently, $a=b=0$.

Therefore, if $\tau=(-1,1,1,-1)$ or $\tau=(1,-1,-1,1)$, then $A_{\tau}$ does not have automorphisms of degree $\geq 3$. Consequently,
$$\varphi=(\alpha_1x_1+\beta_1x_2+\gamma_1,\,\alpha_2x_1+\beta_2x_2+\gamma_2 ),\,\,\alpha_i,\,\beta_i,\,\gamma_i \in K.$$
By (\ref{f12}), we have
 $$((x_i\otimes x_j)\tau)\varphi\otimes \varphi= ((x_i\otimes x_j)\varphi\otimes \varphi)\tau,
\, i, j \in\{1,\,2\}.$$
Hence
\bes
q_{ij}(\alpha_jx_1+\beta_jx_2+\gamma_j)\otimes(\alpha_ix_1+\beta_ix_2+\gamma_i)\\
=((\alpha_ix_1+\beta_ix_2+\gamma_i)\otimes
(\alpha_jx_1+\beta_jx_2+\gamma_j))\tau, \ \  i,j \in \{1,2\}.
\ees
Comparing the coefficients in the terms $x_1\otimes x_1, \, x_1\otimes x_2, \, x_2\otimes x_1,\,x_2\otimes x_2,
\,x_1,\,x_2$, and $1$ in this equality, we get
$$(q_{ij}-q_{11})\alpha_i\alpha_j=(q_{ij}-q_{12})\alpha_i\beta_j=(q_{ij}-q_{21})\alpha_j\beta_i=(q_{ij}-q_{22})\beta_i\beta_j$$
$$=(q_{ij}-1)(\alpha_i\gamma_j+\alpha_j\gamma_i)=(q_{ij}-1)(\beta_i\gamma_j+\beta_j\gamma_i)=(q_{ij}-1)\gamma_i\gamma_j=0,\,i, j \in\{1,\,2\}.
$$
Varying the values of $i,j\in \left\{1,2\right\}$, we have
\begin{equation}\label{f18}
\begin{split}
(q_{11}-q_{12})\alpha_1\beta_1=(q_{11}-q_{21})\alpha_1\beta_1=(q_{11}-q_{22})\beta^2_1=(q_{ 11}-1)\alpha_1\gamma_1\\=
(q_{11}-1)\beta_1\gamma_1=(q_{11}-1)\gamma^2_1=0,
\end{split}
\end{equation}

\begin{equation}\label{f19}
\begin{split}
(q_{12}-q_{11})\alpha_1\alpha_2=(q_{12}-q_{21})\alpha_2\beta_1=(q_{12}-q_{22})\beta_1\beta_2\\=(q_{12}-1)(\alpha_1\gamma_2+\alpha_2\gamma_1)=(q_{12}-1)(\beta_1\gamma_2+\beta_2\gamma_1)=(q_{12}-1)\gamma_1\gamma_2=0,
\end{split}
\end{equation}

\begin{equation}\label{f20}
\begin{split}
(q_{21}-q_{11})\alpha_1\alpha_2=(q_{21}-q_{12})\alpha_2\beta_1=(q_{21}-q_{22})\beta_1\beta_2\\= (q_{21}-1)(\alpha_2\gamma_1+\alpha_1\gamma_2)=(q_{21}-1)(\beta_2\gamma_1+\beta_1\gamma_2)=(q_{21}-1)\gamma_2\gamma_1=0,
\end{split}
\end{equation}

\begin{equation}\label{f21}
\begin{split}
(q_{22}-q_{11})\alpha^2_2=(q_{22}-q_{12})\alpha_2\beta_2=(q_{22}-q_{21})\alpha_2\beta_2=(q_{22}-1)\alpha_2\gamma_2\\=
(q_{22}-1)\beta_2\gamma_2=(q_{22}-1)\gamma^2_2=0.
\end{split}
\end{equation}

If $\tau=(-1, 1, 1, -1)$ or $\tau=(1, -1, -1, 1)$, it follows from (\ref{f18}) -- (\ref{f21}) that
\bes
\alpha_1\beta_1=\alpha_1\alpha_2=\beta_1\beta_2=\alpha_2\beta_2=\gamma^2_1=\gamma^2_2=0
\ees
 or
 \bes
 \alpha_1\beta_1=\alpha_1\alpha_2=\beta_1\beta_2=\alpha_2\beta_2=\alpha_1\gamma_2+\alpha_2\gamma_1=\beta_1\gamma_2+\beta_2\gamma_1=\gamma_1\gamma_2=0,
 \ees
 respectively.
Using this and (\ref{f6}), we get
\bes
\alpha_1\neq 0, \, \beta_2\neq 0, \, \beta_1=\alpha_2=\gamma_1=\gamma_2=0
\ees
 or
 \bes
 \beta_1\neq 0, \, \alpha_2\neq 0, \,  \alpha_1= \beta_2=\gamma_1=\gamma_2=0.
  \ees
  This means that  $$\varphi=(\alpha_1x_1, \ \beta_2x_2) \,\,  \text{or} \, \,  \varphi=(\beta_1x_2, \ \alpha_2x_1).$$
Using Lemma \ref{1}, it is easy to check that each of these automorphisms is an automorphism of $A_\tau$. Consequently, $\mathrm {Aut} A_\tau = G_3$.  $\Box$

\begin{lm}\label{6}
If $\tau=(-1,-1,-1,-1)$ then $\mathrm{Aut} A_\tau = G_2$.
\end{lm}
\Proof Let $u, v \in X^*$, $s=\mathrm{d}(u)$, and $t=\mathrm{d}(u)$. If $\tau=(-1,-1,-1,-1)$ then, by Lemma 1,
\bes
(u\otimes v)\tau=(-1)^{st}(v\otimes u).
\ees
Notice that $(-1)^{st}=-1$ if and only if both $s$ and $t$ are odd.
This implies that
\bee\label{f22}
(f\otimes g)\tau=(-1)^{ij}(g\otimes f)
\eee
for any homogeneous elements $f\in B_i$ and $g\in B_j$ with respect to the grading (\ref{u1}).

Let $\varphi\in G_2=G_{odd}$ be an arbitrary odd automorphism of $A=K\langle x_1,x_2\rangle$. Notice that $\varphi(B_i)\subseteq B_i$. Using (\ref{f22}),
we get
\bes
(f\otimes g)\tau(\varphi\otimes\varphi)=(-1)^{ij}(g\otimes f)(\varphi\otimes\varphi)\\
=(-1)^{ij}(\varphi(g)\otimes \varphi(f))=(\varphi(f)\otimes \varphi(g))\tau=
(f\otimes g)(\varphi\otimes\varphi)\tau,
\ees
that is, (\ref{f12}) holds for $\varphi$. Consequently, $\varphi\in \mathrm{Aut} A_\tau$ and $G_2\subseteq \mathrm{Aut} A_\tau$.

Let $\varphi=(f_1,\, f_2)$ be an arbitrary automorphism of the algebra $A_\tau$. We prove that $\varphi\in G_2$ by induction on $\deg(\varphi)$.
Let $\deg(\varphi)=2$,  that is,
$$\varphi=(\alpha_1x_1+\beta_1x_2+\gamma_1,\alpha_2x_1+\beta_2x_2+\gamma_2),  \ \  \alpha_i,\beta_i,\gamma_i \in K.$$
It follows from (\ref{f18}) -- (\ref{f21}) that
$\gamma_1=\gamma_2=0$. Consequently,

$$\varphi=(\alpha_1x_1+\beta_1x_2,\alpha_2x_1+\beta_2x_2)\in G_2.$$

Suppose that $\deg(\varphi)\geq 3$, $\deg f_1=s_1$,  and $\deg f_2=s_2$.
Then  (\ref{f14}) -- (\ref{f17}) imply that
$$(-1-(-1)^{s_is_j})a^{s_i+s_j}=(-1-(-1)^{s_is_j})a^{s_i}b^{s_j}$$
$$=(-1-(-1)^{s_is_j})a^{s_j}b^{s_i}=(-1-(-1)^{s_is_j})b^{s_i+s_j}=0, \, i,\,j\in\{1,\,2\}.$$
If $s_1s_2$ is even then these equations give $a=b=0$ over a field of characteristic $\neq 2$. Consequently, $s_1$ and $s_2$ are both odd.
Suppose that $s_1|s_2$ for definiteness by Lemma \ref{2}. Set $t=s_2/s_1$. Then there exists an odd automorphism $\lambda=(x_1,x_2+\b x_1^t)$ such that
$\deg(\varphi\lambda)<\deg(\varphi)$. Since $\lambda\in G_2\subseteq A_\tau$, it follows that $\varphi\lambda\in A_\tau$. By the induction proposition,
$\varphi\lambda\in G_2$. Consequently, $\varphi\in G_2$. $\Box$

\begin{lm}\label{8}
If $\tau=(1,1,1,-1)$ then $\mathrm {Aut}_K A_\tau = G_4$.
\end{lm}
\Proof Consider the grading
\bes
A=C_0\oplus C_1
\ees
of the algebra $A$, where $C_0$ is the linear span of all monomials of even degree in $x_2$ and $C_1$ is the linear span of all monomials of odd degree in $x_2$.

Let $u,v \in X^*$, $\mathrm{mdeg}(u)=(s,s')$, $\mathrm{mdeg}(v)=(t,t')$. Lemma 1 gives
\bes
(u\otimes v)\tau=(-1)^{s't'}(v\otimes u).
\ees
Notice that $(-1)^{s't'}=-1$ if and only if both $s$ and $t$ are odd.
This implies that
\bee\label{f23}
(f\otimes g)\tau=(-1)^{ij}(g\otimes f)
\eee
for any homogeneous elements $f\in C_i$ and $g\in C_j$.

Let $\varphi\in G_4$. Notice that $\varphi(C_i)\subseteq C_i$. Using (\ref{f23}), we easily get
\bes
(f\otimes g)\tau(\varphi\otimes\varphi)=(-1)^{ij}(g\otimes f)(\varphi\otimes\varphi)\\
=(-1)^{ij}(\varphi(g)\otimes \varphi(f))=(\varphi(f)\otimes \varphi(g))\tau=
(f\otimes g)(\varphi\otimes\varphi)\tau.
\ees
 Consequently, $\varphi\in \mathrm{Aut} A_\tau$ and $G_4\subseteq \mathrm{Aut} A_\tau$.

Let $\varphi=(f_1,f_2)$ be an arbitrary automorphism of the algebra $A_\tau$. We prove by induction on $\deg\varphi$ that $\varphi\in G_4$.
Let $\deg(\varphi)=2$, i.e.,
$$\varphi=(\alpha_1x_1+\beta_1x_2+\gamma_1,  \alpha_2x_1+\beta_2x_2+\gamma_2), \ \ \ \alpha_i, \beta_i, \gamma_i \in K.$$

It follows from (\ref{f18}) -- (\ref{f21}) that $$\beta^2_1=\beta_1\beta_2=\alpha^2_2=\alpha_2\beta_2=\alpha_2\gamma_2=\beta_2\gamma_2=\gamma^2_2=0.$$
Using this and (\ref{f6}), we get
$$\alpha_1\neq 0, \, \beta_2\neq 0, \, \beta_1=\alpha_2=\gamma_2=0.$$
This means that
$$\varphi=(\alpha_1x_1+\gamma_1,\,\beta_2x_2)\in G_4.$$

 Suppose that  $\deg(\varphi)\geq 3$, $\deg f_1=s_1$, and $\deg f_2=s_2$.
As in the proof of Lemma 4, (\ref{f14}) -- (\ref{f17}) imply that
$$(-1-1^{s^2_2})a^{2s_2}=(1-(-1)^{s^2_1})b^{2s_1}=(1-(-1)^{s_1s_2})b^{s_1+s_2}=(-1-(-1)^{s^2_2})b^{2s_2}=0.$$
From this we immediately get $a=0$ over a field of characteristic $\neq 2$. Then $b\neq 0$. This is possible only if $(s_1,s_2)=(2k,2l-1)$ for some positive integers $k,l$. By Lemma \ref{2}, $s_2|s_1$ and $t=s_1/s_2=2k/(2l-1)$ is even. Consequently, there exists
$\lambda=(x_1+\b x_2^t,x_2)\in G_4$ such that
$\deg(\varphi\lambda)<\deg(\varphi)$. Since $\lambda\in G_4\subseteq A_\tau$, it follows that $\varphi\lambda\in A_\tau$. By the induction proposition,
$\varphi\lambda\in G_4$. Consequently, $\varphi\in G_4$. $\Box$

\begin{lm}\label{9}
If $\tau=(1,-1,-1,-1)$ then $\mathrm{Aut} A_\tau = G_5$.
\end{lm}
\Proof
 Let $\varphi=(f_1,\, f_2)$ be an automorphism of the algebra $A_\tau$. If $\deg(\varphi)\geq 3$, $\deg f_1=s_1$, and $\deg f_2=s_2$,
 then (\ref{f14}) -- (\ref{f17}) imply that
$$(-1-1^{s_1s_2})a^{s_1+s_2}=(1-(-1)^{s^2_1})b^{2s_1}$$
$$=(-1-(-1)^{s_1s_2})b^{s_1+s_2}=(-1-(-1)^{s^2_2})b^{2s_2}=0.$$
It is easy to check that over a field of characteristic $\neq 2$, these equations give $a=b=0$. Therefore, $A_\tau$ does not have 
any automorphisms with degree $\geq 3$.

Suppose that
$$\varphi=(\alpha_1x_1+\beta_1x_2+\gamma_1,\,\alpha_2x_1+\beta_2x_2+\gamma_2 ),\,\,\alpha_i,\,\beta_i,\,\gamma_i \in K.$$
It follows from (\ref{f18}) -- (\ref{f21}) that
$$\alpha_1\beta_1=\beta^2_1=\alpha_1\alpha_2=\alpha_1\gamma_2+\alpha_2\gamma_1=\beta_1\gamma_2+\beta_2\gamma_1=\gamma_1\gamma_2=\alpha^2_2=\alpha_2\gamma_2=\beta_2\gamma_2=\gamma^2_2=0.$$
Using this and (\ref{f6}), we get
$$\alpha_1\neq0, \, \beta_2\neq0, \, \beta_1=\alpha_2=\gamma_1=\gamma_2=0.$$
 This means that $$\varphi=(\alpha_1x_1,\, \beta_2x_2)\in G_5.$$
Using Lemma \ref{1}, it is easy to check that each element of $G_5$ is an automorphism of $A_\tau$. Consequently, $\mathrm {Aut} A_\tau=G_5$. $\Box$

\begin{lm}\label{10}
If $\tau=(q_{11},q_{12},q_{21},q_{22})$ with $q_{12}\neq\pm1$ then $\mathrm{Aut} A_\tau = G_{5}$.
\end{lm}
\Proof
Let $\tau=(q_{11},q_{12},q_{21},q_{22})$ with $q_{12}\neq\pm 1$.
We have $q^2_{11}=q^2_{22}=q_{12}q_{21}=1$ since $\tau^2={\rm id}$. Consequently,  $q_{11}=\pm1,\, q_{22}=\pm1,\, q_{21}\neq q_{12},\, q_{12}\neq\pm1$.

Let $\varphi=(f_1,\, f_2)$ be an automorphism of the algebra $A_\tau$. If $\deg(\varphi)\geq 3$, $\deg f_1=s_1$, and $\deg f_2=s_2$,
 then   over a field of characteristic $\neq 2$,  (\ref{f14}) -- (\ref{f17}) give $a=b=0$. Therefore, $A_\tau$ does not have any automorphisms with degree $\geq 3$.

Suppose that
$$\varphi=(\alpha_1x_1+\beta_1x_2+\gamma_1,\,\alpha_2x_1+\beta_2x_2+\gamma_2 ),\,\,\alpha_i,\,\beta_i,\,\gamma_i \in K.$$
It follows from  (\ref{f18}) -- (\ref{f21}) that
 $$\alpha_1\beta_1=\alpha_1\alpha_2=\alpha_2\beta_1=\beta_1\beta_2=\alpha_2\beta_2=\alpha_2\gamma_1+\alpha_1\gamma_2=\beta_2\gamma_1+\beta_1\gamma_2=\gamma_1\gamma_2=0.$$
Using this and (\ref{f6}), we get
$$\alpha_1\neq0, \, \beta_2\neq0, \, \beta_1=\alpha_2=\gamma_1=\gamma_2=0.$$
 This means that $$\varphi=(\alpha_1x_1,\, \beta_2x_2)\in G_5.$$
Using Lemma \ref{1}, it is easy to check that each element of $G_5$ is an automorphism of $A_\tau$. Consequently, $\mathrm {Aut} A_\tau=G_5$. $\Box$

\begin{theor}\label{t1}
Let $A_{\tau}=\langle K\langle x_1,\, x_2\rangle, \tau\rangle$ be a free braided associative algebra in two generators $x_1,x_2$ over a field $K$ of arbitrary characteristic $\neq2$ with an involutive diagonal braiding
$\tau=(q_{11},q_{12},q_{21},q_{22})$.  Then

$(1)$ $\mathrm{Aut}A_\tau = \mathrm{Aut}K\langle x_1,\, x_2\rangle$ if $q_{ij}=1$ for all $i,j$;

$(2)$ $\mathrm{Aut}A_\tau = G_{odd}$  if $q_{ij}=-1$ for all $i,j$, where $G_{odd}$ is the group of all odd automorphisms of $K\langle x_1,\, x_2\rangle$;

$(3)$ $\mathrm{Aut}A_\tau \cong (K^*\times K^*) \rtimes \mathbb{Z}_2$ if $q_{11}=q_{22}$, $q_{12}=q_{21}$, and $q_{11}q_{12}=-1$;

$(4)$ $\mathrm{Aut}A_\tau \cong G_4$ if $q_{12}=1$ and $q_{11}q_{22}=-1$;

$(5)$ $\mathrm{Aut}A_\tau \cong K^*\times K^*$ if $q_{12}\neq \pm 1$ or $q_{12}=-1, q_{11}q_{22}=-1$.
\end{theor}
\Proof
The statements $(1)$ and $(2)$ are proven in Lemmas \ref{3} and \ref{6}, respectively. Let $q_{11}=q_{22}$, $q_{12}=q_{21}$, and $q_{11}q_{12}=-1$. Using Proposition  \ref{p1}, we can assume that $\tau=(-1,1,1,-1)$ or $\tau=(1,-1,-1,1)$. By Lemma \ref{4}, $\mathrm{Aut}A_\tau \cong G_3$. Notice that
\bes
G_3\cong (K^*\times K^*) \rtimes \mathbb{Z}_2,
\ees
where $\mathbb{Z}_2$ is the subgroup of $G_3$ generated by $(x_2,x_1)$.

If $q_{12}=1$ and $q_{11}q_{22}=-1$ then we can assume that $\tau=(1,1,1,-1)$ by Proposition \ref{p1}. By Lemma \ref{8}, $\mathrm{Aut}A_\tau \cong G_4$.
If $q_{12}\neq \pm 1$ then $\mathrm{Aut}A_\tau = G_{5}\cong K^*\times K^*$ by Lemma \ref{10}. If $q_{12}=-1$ and $q_{11}q_{22}=-1$ then we can assume that
$\tau=(1,-1,-1,-1)$ by Proposition \ref{p1}. Then Lemma \ref{9} gives $\mathrm{Aut}A_\tau = G_{5}$ again. $\Box$

Since $\mathrm{Aut}K\langle x_1,\, x_2\rangle\cong \mathrm{Aut}K[x_1,x_2]$, the groups $\mathrm{Aut}K\langle x_1,\, x_2\rangle$, $K^*\times K^*$, and 
$(K^*\times K^*) \rtimes \mathbb{Z}_2$ are the automorphism groups of quantum planes \cite{AC}.

\bigskip

{\bf Acknowledgments}

\bigskip

This work was supported by the grant of the Ministry of Education and Science of the Republic of Kazakhstan (project AP08052290).

\end{document}